\documentclass[12pt,twoside]{amsart}
\usepackage{amssymb}
\usepackage{latexsym}
\usepackage{amsfonts}
\usepackage{amsmath}
\usepackage{amssymb}

\newcommand{\K}{\mathcal K}

\newcommand{\E}{\mathcal E}

\newtheorem{theorem}{Theorem}[section]
\newtheorem{proposition}[theorem]{Proposition}

\theoremstyle{definition}

\newtheorem{definition}[theorem]{Definition}
\newtheorem{example}[theorem]{Example}

\newcommand{\sy}[1]{{\sf S}_{#1}}

\newcommand{\sym}[1]{{\sf Sym}\,#1}
\renewcommand{\wr}{\,{\sf wr}\,}

\renewcommand{\leq}{\leqslant}
\renewcommand{\geq}{\geqslant}

\newcommand{\aut}[1]{{\sf Aut}(#1)}

\newcommand{\psl}[2]{{\sf PSL}_{#1}(#2)}

\newcommand{\pgammal}[2]{{\sf P}\Gamma{\sf L}_{#1}(#2)}

\newcommand{\alt}[1]{{\sf A}_{#1}}
\newcommand{\Alt}[1]{{\sf Alt}\,#1}
\newcommand{\dih}[1]{{\sf D}_{#1}}

\begin{document}

\title[Cartesian decompositions for permutation groups]{Identifying
Cartesian decompositions preserved by transitive permutation groups}
\author{Robert W. Baddeley, Cheryl E. Praeger and Csaba Schneider}
\address[Baddeley]{32 Arbury Road,
  Cambridge CB4 2JE, UK}
                                                                                              
\address[Praeger \& Schneider]{School of Mathematics and Statistics\\
The University of Western Australia\\
35 Stirling Highway, Crawley WA 6009\\
Australia}
\email{robert.baddeley@ntlworld.com, praeger@maths.uwa.edu.au,\newline
csaba@maths.uwa.edu.au\protect{\newline} {\it WWW:}
www.maths.uwa.edu.au/$\sim$praeger, www.maths.uwa.edu.au/$\sim$csaba}


\date{\today}
\thanks{This paper forms part of an Australian Research Council large
grant project.}

\maketitle

\section{Introduction}

Intuitively, a Cartesian decomposition of a finite set $\Omega$ is a way of
identifying $\Omega$ with a Cartesian product
$\Gamma_1\times\dots\times\Gamma_\ell$ of smaller sets $\Gamma_i$.  However we
do not wish to distinguish between two such identifications if the second can be
obtained from the first by re-naming the elements in the individual sets
$\Gamma_i$.  Nor do we wish to distinguish between, say,
$\Gamma_1\times\Gamma_2$ and $\Gamma_2\times\Gamma_1$.  Thus by a Cartesian
decomposition we will mean an equivalence class of identifications of $\Omega$
with a Cartesian product  under a certain
notion of equivalence.  Our formal Definition~\ref{csdef} encompasses these
ideas, but at first reading this may not be apparent.  We therefore develop the
concept further in Section~\ref{sect:cd} before giving the formal definition.

Our aim is to describe the theory of Cartesian decompositions preserved by some
member of a large family of finite transitive permutation groups called innately
transitive groups.  Innately transitive groups are defined in
Section~\ref{sect:it}, and for such a group $G$, the Cartesian decompositions
preserved by $G$ correspond to certain families of subgroups, called Cartesian
systems, of a normal subgroup $M$ of $G$.  Many Cartesian decompositions
correspond to direct decompositions of $M$.  These are called $M$-normal and are
defined in Section~\ref{sect:ncd}.  The non-normal $G$-invariant Cartesian
decompositions occur for $M$ of the form $M=T^k$, where $T$ is a nonabelian
simple group, and are often related to factorisations of $T$.  The various types
of such non-normal Cartesian decompositions are discussed in
Section~\ref{sect:factns}, and simple examples illustrating most of the
possibilities are given in Section~\ref{sect:ex}.

\section{The concept of Cartesian decompositions}\label{sect:cd}

In this section we develop the concept of a Cartesian decomposition,
giving a formal definition in Definition~\ref{csdef}.  An identification of a
finite set $\Omega$ with a Cartesian product
$\Omega_1\times\cdots\times\Omega_{\ell}$ is a bijection
$\varphi:\Omega\rightarrow\Omega_1\times\cdots\times\Omega_{\ell}$.  A second
bijection $\varphi':\Omega\rightarrow\Delta_1\times\cdots\times\Delta_{\ell'}$
is defined to be {\em equivalent} to $\varphi$ if $\ell=\ell'$, and there
exist a permutation $\pi\in\sy\ell$ and bijections
$\alpha_1:\Omega_1\rightarrow\Delta_{1\pi},\ldots,\alpha_\ell:\Omega_\ell\rightarrow\Delta_{\ell\pi}$
such that
$$
\varphi'(\omega)=(\alpha_1\times\cdots\times\alpha_\ell)_\pi(\varphi(\omega))\quad\mbox{for
all}\quad\omega\in\Omega,
$$
where $(\alpha_1\times\cdots\times\alpha_\ell)_\pi$ denotes the bijection
$\Omega_1\times\cdots\times\Omega_\ell\rightarrow\Delta_1\times\cdots\times\Delta_\ell$ defined by
$$
(\alpha_1\times\cdots\times\alpha_\ell)_\pi(\omega_1,\ldots,\omega_\ell)=
\left(\alpha_{1\pi^{-1}}(\omega_{1\pi^{-1}}),\ldots,\alpha_{\ell\pi^{-1}}
(\omega_{\ell\pi^{-1}})\right).
$$
Whether or not a second bijection is equivalent to $\varphi$ can be decided 
using the preimages of the natural projection maps
$\sigma_i:\Omega_1\times\cdots\times\Omega_\ell\rightarrow\Omega_i$ as follows.
For $i=1,\ldots,\ell$ set
$$
\Gamma_i=\left\{\{\omega\ |\ \sigma_i(\varphi(\omega))=\omega_i\}\ |\
\omega_i\in\Omega_i\right\}.
$$
It is easy to see that each of the $\Gamma_i$ is a partition of
$\Omega$, and that
\begin{equation}\label{cseq}
|\gamma_1\cap\cdots\cap\gamma_\ell|=1\quad\mbox{for all}\quad
 \gamma_1\in\Gamma_1, \ldots,\gamma_\ell\in\Gamma_\ell.
\end{equation}
Moreover two bijections
$\varphi$ and $\varphi'$ are equivalent under the equivalence
relation defined above if and only if they give rise to the
same set of partitions of $\Omega$; the proof of this is left to the reader. Hence each
equivalence class of Cartesian decompositions gives rise to a unique set of
partitions satisfying~\eqref{cseq}. Conversely,
if $\Gamma_1,\ldots,\Gamma_\ell$ are partitions of $\Omega$ such
that~\eqref{cseq} holds, then $\Omega$ can be identified with the
Cartesian product $\Gamma_1\times\cdots\times\Gamma_\ell$ as
follows. Let $\omega\in\Omega$, and let
$\gamma_1\in\Gamma_1, \ldots,\gamma_\ell\in\Gamma_\ell$ such that
$\omega\in\gamma_1\cap\cdots\cap\gamma_\ell$. Such blocks $\gamma_i$
exist and are unique because each $\Gamma_i$ is a partition of $\Omega$. Define
$\psi:\Omega\rightarrow\Gamma_1\times\cdots\times\Gamma_\ell$ by
$\psi(\omega)=(\gamma_1,\ldots,\gamma_\ell)$. Property~\eqref{cseq}
ensures that $\psi$ is a bijection.

Now the motivation behind the following definition should be clear.

\begin{definition}\label{csdef}
If $\Omega$ is a finite set, then a set
$\E=\{\Gamma_1,\ldots,\Gamma_\ell\}$ of partitions of $\Omega$ is said to
be a
{\em Cartesian decomposition} of $\Omega$ if~\eqref{cseq} holds.
\end{definition}

This definition of Cartesian decompositions enables us to study
Cartesian decompositions that are invariant under the action of a
permutation group. If $\E$ is a Cartesian decomposition of $\Omega$
and $g\in\sym\Omega$ then we say that $\E$ is invariant under $g$ if
the partitions in $\E$ are permuted by $g$. The stabiliser in
$\sym\Omega$ of a
Cartesian decomposition $\E$ is obviously a subgroup of $\sym\Omega$.

We give two simple examples of Cartesian decompositions that are invariant
under the action of some transitive permutation group.

\begin{example}\label{ex}
Let $\Omega=\{1,2\}\times \{1,2,3\}$. 
Then the identity map $(i,j)\mapsto (i,j)$ of $\Omega$ is a
bijection whose corresponding Cartesian
decomposition contains the two partitions given by the columns and
the rows of the following grid.
\begin{center}
$$\begin{array}{|c|c|c|}
\hline
(1,1) & (1,2) & (1,3) \\
\hline
(2,1) & (2,2) & (2,3)\\
\hline
\end{array}
$$
\end{center}
Hence the Cartesian decomposition corresponding to the identity
map on $\Omega$ consists of the following two partitions, namely the rows of the
grid,
$$
\Gamma_1=\left\{\{(1,1),(1,2),(1,3)\},\{(2,1),(2,2),(2,3)\}\right\}
$$
and the columns of the grid
$$
\Gamma_2=\left\{\{(1,1),(2,1)\},\{(1,2),(2,2)\},\{(1,3),(2,3)\}\right\}.
$$
Note that the two partitions in this Cartesian decomposition have
different sizes. Such a Cartesian decomposition is said to be {\em inhomogeneous}.
The group $G=\sy 2\times\sy 3$ in its natural action on $\Omega$ is the full
stabiliser of $\{\Gamma_1,\Gamma_2\}$ in $\sym\Omega$, and is transitive on
$\Omega$, but, as this Cartesian decomposition is inhomogeneous, 
no element of $G$ swaps $\Gamma_1$ and $\Gamma_2$. Hence
$G$ is intransitive on the Cartesian decomposition
$\{\Gamma_1,\Gamma_2\}$ of $\Omega$.
\end{example}

\begin{example}\label{exhom}
Let $\Gamma$ be a finite set, let $\ell\geq 2$, and let $\Omega=\Gamma
\times\dots\times\Gamma=\Gamma^\ell$. The wreath product
$W=\sym\Gamma\wr\sy\ell$ is the semidirect product of its normal 
subgroup $N=(\sym\Gamma)^\ell$ and a subgroup $H\cong S_\ell$. The product
action of $W$ on $\Omega$ is defined by
$$
(\gamma_1,\dots,\gamma_\ell)^{xh}=(\gamma_{1h^{-1}}^{x_{1h^{-1}}},\dots,
\gamma_{\ell h^{-1}}^{x_{\ell h^{-1}}})
$$
for all $(\gamma_1,\dots,\gamma_\ell)\in\Omega$, $x=(x_1,\dots,
x_\ell)\in N$, and $h\in H$, where we write the image of $\gamma
\in\Gamma$ under $y\in\sym\Gamma$ as $\gamma^y$. Clearly $W$ is transitive
on $\Omega$.
The Cartesian decomposition corresponding to the identity map on $\Omega$
is $\E=\{\Gamma_1,\ldots,\Gamma_\ell\}$
,
where  $\Gamma_i$ is the partition of $\Omega$ into disjoint
subsets according to the $i$-th coordinate of a point in
$\Omega$, that is to say, the parts of $\Gamma_i$ are indexed by
$\Gamma$ and the $\gamma$-part is the set of all points $(\gamma_1,\dots,\gamma_\ell)$
with $\gamma_i=\gamma$. Thus $|\Gamma_i|=|\Gamma|$ for all
$i$. A Cartesian decomposition
$\{\Gamma_1,\ldots,\Gamma_\ell\}$ for which the $\Gamma_i$ all have the
same cardinality is said to be {\em homogeneous}. Thus $\E$ is homogeneous.
Also each element $xh\in W$ maps the partition $\Gamma_i$ to the partition
$\Gamma_{ih}$. Thus $W$ preserves the Cartesian decomposition $\E$.
In fact $W$ is the full stabiliser of $\E$ in $\sym\Omega$ and $W$ permutes
the partitions $\Gamma_i$ transitively.
\end{example}

A Cartesian decomposition $\E$ is called \emph{$G$-transitive} if it is
$G$-invariant and $G$ acts transitively on the set $\E$ of partitions.  Thus the
Cartesian decomposition of Example~\ref{exhom} is $W$-transitive but not
$N$-transitive.

\section{Innately transitive groups}\label{sect:it}

The class of primitive permutation groups that preserve a Cartesian
decomposition of the underlying set is well-understood, and is
described in~\cite{prae:inc}. The original aim of the research presented in this article was to
extend this result to describe Cartesian
decompositions that are preserved by a permutation group in which all
minimal normal subgroups are transitive. Such a group is said to be
{\em quasiprimitive}. We found that the methods used to
achieve this goal give a description of Cartesian
decompositions preserved by a larger class of groups, namely those
that have at least one transitive minimal normal subgroup. Such a
group is said to be {\em innately transitive} and a transitive minimal normal
subgroup of an innately transitive group is called a {\em
plinth}. In particular each primitive or quasiprimitive group is
innately transitive, and the class of innately transitive groups also
contains many interesting non-quasiprimitive groups.
Innately transitive groups are studied in~\cite{bamprae}.
The essential reason why they provide us with the right
context for our research is the following result.

\begin{proposition}\label{prop}
Let $G$ be an innately transitive group acting on $\Omega$, and let
$M$ be a plinth of $G$. If $\{\Gamma_1,\ldots,\Gamma_\ell\}$ is a
$G$-invariant Cartesian decomposition of $\Omega$, then each
$\Gamma_i$ is an $M$-invariant partition of $\Omega$.
\end{proposition}

A plinth $M$ of an innately transitive group $G$ on $\Omega$ is
transitive on $\Omega$. It is well-known (see, for example,
\cite[Theorem 1.5A]{DM96}) that, for a transitive permutation group
$M$ on $\Omega$ and a fixed $\omega\in\Omega$, there is a one-to-one
correspondence between the set of $M$-invariant partitions of $\Omega$
and the set of subgroups $K$ of $M$ containing the stabiliser $M_\omega$. 
Thus the subgroups $K_1,\ldots,K_\ell$ of $M$ corresponding to 
$M$-invariant partitions $\Gamma_1,\ldots,\Gamma_\ell$ of $\Omega$
are sufficient to determine the partitions $\Gamma_i$,
and we can decide from certain properties of the $K_i$ whether 
or not the $\Gamma_i$ form a Cartesian decomposition of $\Omega$. 

More precisely, let $G$ be an innately transitive group acting on
$\Omega$ with plinth $M$, and let $\{\Gamma_1,\ldots,\Gamma_\ell\}$ be
a $G$-invariant Cartesian decomposition of $\Omega$. By
Proposition~\ref{prop}, each of the
$\Gamma_i$ is an $M$-invariant partition of $\Omega$. Fix
$\omega\in\Omega$, let
$\gamma_1\in\Gamma_1,\ldots,\gamma_\ell\in\Gamma_\ell$ be such that
$\{\omega\}=\gamma_1\cap\cdots\cap\gamma_\ell$, and, for
$i=1,\ldots,\ell$, set $K_i=M_{\gamma_i}$. It is proved 
in~\cite[Lemma 2.2]{recog}
that the set $\{K_1,\ldots,K_\ell\}$ is invariant under conjugation by
$G_\omega$ and has the following two important properties:
\begin{eqnarray}
\label{cs1}\bigcap_{i=1}^\ell K_i&=&M_\omega,\\
\label{cs2}K_i\left(\prod_{j\neq i}K_j\right)&=&M\quad\mbox{for all }i\in\{1,\ldots,\ell\}.
\end{eqnarray}

\begin{definition}\label{def3.2}
Let $G$ be a transitive permutation group on $\Omega$ with plinth $M$ and let
$\K=\{K_1,\ldots,K_\ell\}$ be a $G_\omega$-invariant set of subgroups
of $M$ such that~\eqref{cs1}
and~\eqref{cs2} hold. Then $\K$ is said to be a {\em Cartesian system} of
subgroups in $M$ with respect to $\omega$. 
\end{definition}

Conversely, using the correspondence explained above, any Cartesian system in
$M$ leads to a $G$-invariant Cartesian decomposition.  Thus the set of Cartesian
decompositions invariant under the action of an innately transitive group can be
studied via the set of Cartesian systems in the plinth.

\begin{theorem}\label{one-to-one}
Let $G$ be an innately transitive group acting on a set $\Omega$ with
plinth $M$ and let $\omega$ be a fixed element of $\Omega$. Then there
is a one-to-one correspondence between the set of $G$-invariant
Cartesian decompositions of $\Omega$ and the set of Cartesian systems
in $M$ with respect to $\omega$.
\end{theorem}

Consider Examples~\ref{ex} and \ref{exhom} in terms of Cartesian
systems.

\begin{example}
If $G$ and $\Omega$ are as in Example~\ref{ex}, then $G$ has no
transitive minimal normal subgroup, so $G$ is not innately
transitive. On the other hand, the group $W=NH$ in Example~\ref{exhom}
is innately transitive on $\Omega$, and $W$ preserves 
the Cartesian decomposition $\E=\{\Gamma_1,\ldots,\Gamma_\ell\}$ of
$\Omega$. Suppose that $|\Gamma_i|\geq5$ so that the plinth is 
$M=M_1\times\dots\times M_\ell=(\Alt\Gamma)^\ell$.
Set $\omega=(\gamma,\ldots,\gamma)$. An easy computation shows that 
the Cartesian system corresponding to $\E$ with respect to $\omega$ is
$\{K_1,\ldots,K_\ell\}$ where each $K_i=
(M_i)_\gamma\times \prod_{j\ne i}M_{j}.
$
\end{example}

If $G$ is innately transitive with an abelian plinth then $G$ is primitive (see
\cite{bamprae}) and so, as we mentioned above, all $G$-invariant Cartesian
decompositions have been determined in \cite{prae:inc}.  Thus for the rest of
the paper we will assume that each plinth is nonabelian.

\section{Normal and non-normal Cartesian systems}\label{sect:ncd}

The Cartesian system $\E=\{\Gamma_1,\dots,\Gamma_\ell\}$ presented in
Example~\ref{exhom} has the property that the group $M=(\Alt\Gamma)^\ell$ can be
written as a direct product of $\ell$ subgroups with the action of the $i^{th}$
direct factor corresponding to the $M$-action on the $i^{th}$ partition
$\Gamma_i$.  This is a very useful property, and is perhaps a property possessed
by most of the transitive Cartesian decompositions that might come readily to
mind. We formalise this property of $G$-invariant Cartesian decompositions
for innately transitive groups $G$ with nonabelian plinths.

\begin{definition}\label{def:normal}
Let $G$ be an innately transitive group acting on $\Omega$ with a
non-abelian plinth $M$, and
let $\K=\{K_1,\ldots,K_\ell\}$ be a Cartesian system of subgroups in
$M$ with respect to some $\omega\in\Omega$.
Then $\K$ is said to be {\em $M$-normal} if there are normal subgroups $M_1,\ldots,M_\ell$ of $M$ such
that $M=M_1\times\cdots\times M_\ell$
and each $K_i=(M_i\cap M_\omega)\times\prod_{j\ne i}M_j$.
A $G$-invariant Cartesian decomposition of $\Omega$ is said to be {\em
normal} if the corresponding Cartesian system is $M$-normal for some plinth $M$.
\end{definition}

Normal Cartesian decompositions are considered natural, and they can
be determined using the direct factorisations of the plinth. On the
other hand, not every
Cartesian decomposition is normal. The simplest non-normal Cartesian
decomposition preserved by an innately transitive group is given in
the following example.

\begin{example} Let $G\cong \pgammal 29$ and consider the unique transitive
action of $G$ on a set $\Omega$ of 36~points.  The group $G$ is innately
transitive on $\Omega$, because $G$ has a unique minimal normal subgroup
$M\cong\alt 6$ and $M$ is transitive on $\Omega$.  Moreover if
$\omega\in\Omega$, then $M_\omega\cong\dih{10}$, and it is easy to see that $M$
has subgroups $K_1,\ K_2$, both isomorphic to $\alt 5$, such that $\{K_1,\
K_2\}$ is a Cartesian system of subgroups in $M$ with respect to $\omega$.
Hence $G$ preserves a Cartesian decomposition $\{\Gamma_1,\Gamma_2\}$ of
$\Omega$, where each $\Gamma_i$ has six parts of size $6$.  Since $M$ is simple
and is the unique minimal normal subgroup of $G$, this Cartesian decomposition
cannot be normal.
\end{example}

\section{Factorisations of finite simple groups}\label{sect:factns}

Let $G$ be an innately transitive group on $\Omega$ with a nonabelian plinth
$M$.  If $M$ is simple then it is very unusual for $G$ to preserve a nontrivial
Cartesian decomposition of $\Omega$, any such decomposition is non-normal,
and in fact all such possibilities have been classified explicitly, see
\cite[Theorem 6.1]{recog}.  In this section we outline a theory of non-normal
Cartesian decompositions preserved by innately transitive groups with a nonabelian
plinth, in particular pointing out the role of simple group factorisations.

The group $M$ is a minimal normal subgroup of $G$, and so $M$ is a nonabelian
characteristically simple group and hence is of the form
$M=T_1\times\cdots\times T_k$, where the $T_i$ are finite simple groups, each
isomorphic to the same simple group $T$.  Moreover, the group $G$ acts
transitively by conjugation on the set $\{T_1,\dots,T_k\}$.
For $i=1,\dots,k$, let $\sigma_i:M\rightarrow T_i$
denote the $i$-th projection map.

By Theorem~\ref{one-to-one}, there is a one-to-one correspondence
between the set of $G$-invariant Cartesian decompositions of $\Omega$
and the set of Cartesian systems in $M$ relative to a given point $\omega$.
Since in a Cartesian system $\{K_1,\dots,K_\ell\}$ the
factorisation property~\eqref{cs2} holds, we obtain factorisations of the
simple direct factors of $M$ using the natural projection maps $\sigma_i$, as
follows.  For all $i$, \eqref{cs2}~gives the following
factorisations of $T_i$:
\begin{equation}\label{ti}
T_i=\sigma_i(K_j)\left(\bigcap_{m\neq j}\sigma_i(K_m)\right),\quad\mbox{for all}
\ j=1,\dots,\ell.
\end{equation}
Many of the subgroups $\sigma_i(K_j)$ may coincide with $T_i$,
so we are really interested in the following sets:
\begin{equation}\label{fi}
\mathcal F_i=\{\sigma_i(K_j)\ |\ \sigma_i(K_j)\neq T_i,\
j=1,\ldots,\ell\}.
\end{equation}
The set $\mathcal F_i$ is essentially
independent of $i$, in the sense that if $i_1,\
i_2\in\{1,\ldots,k\}$ then there is some $g\in G$ such that
$T_{i_1}^g=T_{i_2}$, and then we have $\mathcal F_{i_1}^g=\mathcal F_{i_2}$.
In particular $|\mathcal F_i|$ is independent of $i$.

Since the size of $\mathcal F_i$ is an invariant of the corresponding
Cartesian decomposition, one natural way of subdividing the class of
Cartesian decompositions invariant under innately transitive groups
is to use the number $|\mathcal F_i|$. This is achieved in~a forthcoming paper
where Cartesian decompositions in each sub-class are described in
detail. In the rest of this section we summarise the results of~that paper.

Using results on factorisations of finite simple groups in~\cite{bad:fact}
it is easy to prove that $|\mathcal F_i|\leq 3$. Moreover, if
$\mathcal F_i=\{A, B, C\}$ then the information contained in \eqref{ti} is that
\begin{equation}\label{smf}
T_i=A(B\cap C)=B(C\cap A)=C(A\cap B).
\end{equation}
This is called a {\em strong
multiple factorisation} of $T_i$. Such strong
multiple factorisations of finite simple groups are classified
in~\cite{bad:fact}. (All possibilities for $T, A, B, C$ are listed
in~\cite[Table~V]{bad:fact}.) This classification can be used to describe the
$G$-invariant Cartesian decompositions for which the
corresponding $\mathcal F_i$ have 3~elements. If $\mathcal F_i=\{A, B\}$, then
the information contained in \eqref{ti} is precisely that $T=AB$, a
factorisation of the finite
simple group $T$, and moreover each such factorisation may occur
in relation to some Cartesian decomposition; see also
Example~\ref{ex2} in the next section.

If $|\mathcal F_i|=1$, then the corresponding Cartesian decomposition
is either $M$-normal, or the Cartesian system elements contain full diagonal
subgroups isomorphic to $T$ covering exactly two of the simple
direct factors of $M$. If $|\mathcal F_i|=0$, then the $K_i$ are
subdirect subgroups of $M$ and the corresponding Cartesian
decomposition is $M$-normal.

\section{Examples of Cartesian systems}\label{sect:ex}

In this section we give some examples of nontrivial Cartesian systems preserved
by innately transitive groups that illustrate the various sub-classes described
in Section~\ref{sect:factns}.  In these examples $T$ is a finite simple group
and $D(T\times T)$ denotes the straight diagonal subgroup $\{(t,t)\ |\ t\in T\}$
of $T\times T$.  We recall that the sets $\mathcal F_i$
assigned to a Cartesian system are defined in~\eqref{fi}.  The first
example is the smallest one with $\mathcal F_i=\emptyset$, and shows that
every simple group $T$ can occur in this case.

\begin{example}\label{ex1}
Let $G=T\wr\dih{8}=M\rtimes \dih{8}$, where $M=T^4$ and $\dih 8$ acts naturally
on the four~simple direct factors of $M$, that is to say, elements of $\dih 8$
either fix setwise, or interchange, the first two, and the last two, simple 
direct factors of $M$. Let
$$
K_1=D(T\times T)\times T\times T\quad\mbox{and}\quad K_2=T\times T\times
D(T\times T)
$$
so that $K_1\cap K_2=D(T\times T)\times D(T\times T)$ is normalised by $\dih{8}$.
Then the $M$-coset action on
$\Omega=[M:K_1\cap K_2]$ can be extended to $G$ with
$(K_1\cap K_2)\rtimes \dih8$ as the stabiliser of the point $\omega=K_1\cap K_2$.
According to Definition~\ref{def3.2}, $\{K_1,K_2\}$ is
a Cartesian system of subgroups in $M$ with respect to $\omega$.
Hence $G$ preserves the corresponding
Cartesian decomposition of $\Omega$. Clearly in this example we have
$|\mathcal F_i|=0$, and this Cartesian system is $M$-normal. 
\end{example}

In Example~\ref{ex1}, the diagonal subgroups $D(T\times T)$ involved in $K_1$ and
$K_2$ are disjoint in the sense that the diagonal subgroup of $K_1$ is contained
in the direct product of the first two simple direct factors of $M$, while the
diagonal subgroup of $K_2$ is contained in the direct product of the last two
simple direct factors.  It turns out that any two
diagonal subgroups involved in a Cartesian system with $|\mathcal
F_i|=0$ are disjoint.  This also means that each such Cartesian system
is normal.

Next we present two examples with $|\mathcal F_i|=1$, and show that
the class of Cartesian systems with  $|\mathcal F_i|=1$ contains
both normal and non-normal examples. Moreover, the first example shows that,
for each nonabelian simple group $T$ and each of its proper subgroups $A$,
there is such an example with plinth a direct power of $T$ and $\mathcal F_i
=\{A\}$.

\begin{example}
Let $A$ be a proper subgroup of $T$, let  $M=T\times T$,
$G=M\rtimes\sy 2=T\wr\sy 2$,
$K_1=T\times A$, and $K_2=A\times T$. Then the $M$-coset action on
$\Omega=[M:K_1\cap K_2]$ can be extended to $G$ with point stabiliser
$(K_1\cap K_2)\rtimes \sy 2$, and $\{K_1,K_2\}$ is
a Cartesian 
system of subgroups in $M$. Thus $G$ preserves the corresponding
Cartesian decomposition of $\Omega$, 
$|\mathcal F_i|=1$, and the Cartesian system is $M$-normal. Indeed it is
not difficult to see that all Cartesian systems with $|\mathcal
F_i|=1$, and involving no diagonal subgroups, are normal.
\end{example}

\begin{example}
Set $T=\alt 6$,  $A=\alt 5$, and $B=\psl 25$. Then there exists an element 
$\tau\in\aut T$ that swaps $A$ and $B$, such that $\tau^2=1$. 
Let $G_1=M_1\rtimes \sy 2$ where $M_1=T\times T$ and the nontrivial
element $x$ of $\sy 2$ acts via $(t_1,t_2)^x=(t_2^\tau,t_1^\tau)$.
Note that $x$ normalises the subgroup $A\times B$ of $M_1$. Let
$G=G_1\wr\sy 2=(G_1\times G_1)\rtimes\sy 2=M\rtimes \dih{8}$, with $M\cong T^4$
the unique minimal normal subgroup of $G$.
Define 
$$
K_1=A\times B\times D(T\times T)\quad\mbox{and}\quad K_2=D(T\times
T)\times A\times B.
$$
Then the $M$-coset action on
$\Omega=[M:K_1\cap K_2]$ can be extended to $G$ with point stabiliser
$(K_1\cap K_2)\rtimes \dih{8}$, and $\{K_1,K_2\}$ is
a Cartesian 
system of subgroups in $M$. Hence $G$ preserves the corresponding
Cartesian decomposition of $\Omega$, 
$|\mathcal F_i|=1$, but, as the $K_i$ involve diagonal subgroups,
this Cartesian system is not normal. 
\end{example}

Finally we present two further examples, one with $|\mathcal F_i|=2$, and
one with $|\mathcal F_i|=3$.

\begin{example}\label{ex2}
Let $A$ and $B$ be two subgroups of $T$ such that $T=AB$, let $M=T^2$,
$K_1=A\times B$, $K_2=B\times A$, and $G=M\rtimes\sy 2=T\wr\sy 2$.
Then the $M$-coset action on
$\Omega=[M:K_1\cap K_2]$ can be extended to $G$ with point stabiliser
$(K_1\cap K_2)\rtimes\sy 2$, and $\{K_1,K_2\}$ is
a Cartesian 
system of subgroups in $M$. Hence $G$ preserves the corresponding
Cartesian decomposition of $\Omega$,
$|\mathcal F_i|=2$, and the Cartesian system is not normal.
\end{example}

\begin{example}\label{ex3}
Let $T$ be a finite simple group such that $T$, 
$A$, $B$, $C$ form a strong multiple factorisation of $T$, that is to say,
the equations in \eqref{smf} hold. Let $M=T^3$, $G=M\rtimes C_3=T\wr C_3$, 
$$
K_1=A\times B\times C,\quad K_2=B\times C\times A,\quad\mbox{and}\quad
K_3=C\times A\times B.
$$
Then the $M$-coset action on
$\Omega=[M:K_1\cap K_2\cap K_3]$ can be extended to $G$ with point
stabiliser $(K_1\cap K_2\cap K_3)\rtimes C_3$, and $\{K_1,K_2,K_3\}$ is
a Cartesian 
system of subgroups in $M$. Hence $G$ preserves the corresponding
Cartesian decomposition of $\Omega$, 
$|\mathcal F_i|=3$, and the Cartesian system is not normal.
\end{example}

In Example~\ref{ex2}, in addition to the Cartesian system given, there is also
an $M$-normal Cartesian system for the same action of $G$ formed by the
subgroups $(A\cap B)\times T$ and $T\times (A\cap B)$.  Similarly, in
Example~\ref{ex3}, there is an $M$-normal Cartesian system for the same action
of $G$ formed by the subgroups $(A\cap B\cap C)\times T\times T$, $T\times
(A\cap B\cap C)\times T$, and $T\times T\times(A\cap B\cap C)$.  This turns out
to be a rather general phenomenon, and will be studied further in our forthcoming paper.

\end{document}